\documentclass[11pt,a4paper,leqno]{amsart}

\usepackage[latin1]{inputenc}
\usepackage[T1]{fontenc}
\usepackage{amsfonts}
\usepackage{amsmath}
\usepackage{amssymb}
\usepackage{eurosym}
\usepackage{mathrsfs}
\usepackage{palatino}
\usepackage{color}
\usepackage{esint}
\usepackage{url}
\usepackage{verbatim}
\usepackage{graphicx}

\allowdisplaybreaks[4]

\usepackage{enumerate}

\usepackage[pagebackref,hypertexnames=false, colorlinks, citecolor=blue, linkcolor=blue, urlcolor=red]{hyperref}

\newcommand{\R}{\mathbb{R}}
\newcommand{\C}{\mathbb{C}}

\newcommand{\N}{\mathbb{N}}

\newcommand{\Z}{\mathbb{Z}}

\newcommand{\calJ}{\mathcal{J}}

\newcommand{\Dini}{\mathrm{Dini}}
\newcommand{\WBP}{\mathrm{WBP}}
\newcommand{\one}{\mathbf1}

\numberwithin{equation}{section}

\newcommand{\ud}[0]{\,\mathrm{d}}

\newcommand{\norm}[1]{\left\|#1\right\|}

\newcommand{\pair}[2]{\langle #1,#2 \rangle}

\newcommand{\BMO}[0]{\operatorname{BMO}}
\newcommand{\supp}[0]{\operatorname{spt}}

\newcommand{\calD}[0]{\mathcal{D}}

\newcommand{\wt}[1]{{\widetilde{#1}}}
\newcommand{\wh}[1]{{\widehat{#1}}}

\swapnumbers
\theoremstyle{plain}
\newtheorem{theorem}[equation]{Theorem}
\newtheorem{lemma}[equation]{Lemma}
\newtheorem{proposition}[equation]{Proposition}

\theoremstyle{definition}

\theoremstyle{remark}

\title{The $T1$ theorem under Dini regularity}

 \author{Kangwei Li}

 \address[K.L.]{School of Mathematical Sciences, Zhejiang Normal University, Jinhua 321004, China}
 \email{kangwei.li@zjnu.edu.cn}

\makeatletter
\@namedef{subjclassname@2020}{%
  \textup{2020} Mathematics Subject Classification}
\makeatother


\begin{document}
\begin{abstract}
We prove the $T1$ theorem for Calder\'on-Zygmund operators under the classical Dini regularity condition on the kernel, solving a longstanding open problem in singular integral theory. 
\end{abstract}
\maketitle

\section{Introduction}

Given a singular integral operator $T$, a basic question is how to
characterize its $L^2$ boundedness. The classical theory of Calder\'on and
Zygmund \cite{CZ} supplies the analytic framework for this question.
The $T1$ theorem of David and Journ\'e \cite{DJ} gives a particularly
effective answer for nonconvolution operators: for a standard
H\"older-continuous kernel, boundedness on $L^2(\R^d)$ is equivalent
to the weak boundedness property and the conditions
$T1,T^*1\in\BMO(\R^d)$.

We are interested in the amount of regularity needed for this
criterion. Let $K: \R^d \times \R^d \setminus \{(x,x): x\in \R^d\}\to \C$  satisfy 
\begin{equation}\label{eq:size}
 |K(x,y)|\le C_K|x-y|^{-d}
\end{equation}
and
\begin{equation}\label{eq:regularity}
 \begin{split}
 &|K(x+h,y)-K(x,y)|+|K(y,x+h)-K(y,x)|\\
 &\hspace{25mm}\le |x-y|^{-d}\omega\!\left(\frac{|h|}{|x-y|}\right),
 \qquad |h|\le\frac12|x-y|.
 \end{split}
\end{equation}
Here $C_K>0$ and $\omega:[0,1]\to[0,\infty)$ is increasing and subadditive,
with $\omega(0)=0$. We write
\[
 \norm{\omega}_{\Dini_\alpha}
 =\int_0^1\omega(t)\left(1+\log\frac1t\right)^\alpha\frac{dt}{t},
 \qquad \alpha\ge0.
\]
The Dini condition is $\norm{\omega}_{\Dini}<\infty$,
where $\Dini=\Dini_0$. 

The standard H\"older-continuous kernel has $\omega(t)=t^\gamma$ for some $\gamma\in (0,1]$.
Weakening the H\"older hypothesis in the $T1$ theorem has a long
history. The condition $\omega\in\Dini_{1/2}$ is implicit in
Figiel's work \cite{Figiel} and was obtained explicitly by Deng,
Yan and Yang \cite{DYY}. Grau de la Herr\'an and Hyt\"onen
\cite{GH} recovered this threshold through a dyadic representation
and extended it to nonhomogeneous space. Their representation formula
contains model operators with $L^2$ bounds of order the square root
of the complexity. Summing these bounds against
$\omega(2^{-k})$ leads to the condition
\[
 \sum_{k=1}^\infty\sqrt{k}\,\omega(2^{-k})<\infty,
\]
which is precisely the condition $\omega\in\Dini_{1/2}$. 
It is worth mentioning that this idea was further refined in \cite{AMV}, and even extended to multilinear and vector-valued setting. 

We are expecting that Dini regularity condition is sufficient for the $T1$ theorem since 
it matches the Calder\'on--Zygmund theory very well. 
In particular, assuming that the $L^2$ boundedness is known, then the pointwise sparse domination and its weighted
consequences are available at this level of regularity
\cite{Lacey,Lerner}. In \cite[Section~4]{LO}, Lerner and Ombrosi
explicitly raised the question of whether the
$T1$ theorem holds under the Dini condition. The purpose
of this paper is to establish this criterion in the 
Euclidean setting. Thus the regularity required to deduce
boundedness from the $T1$ assumptions agrees with that used in
these subsequent estimates.

Our main result is the following. 

\begin{theorem}\label{thm:main}
Let $T:\mathscr S(\R^d)\to\mathscr S'(\R^d)$ be a continuous linear operator
whose kernel $K$ satisfies \eqref{eq:size} and
\eqref{eq:regularity}. Suppose that
$\omega$ satisfies the Dini condition, 
then $T$ extends to a bounded operator on $L^2(\R^d)$ if and only if
it has the weak boundedness property and
$T1,T^*1\in\BMO(\R^d)$. More precisely,   
\begin{equation}\label{eq:main-bound}
 \norm{T}_{L^2\to L^2}
 \le C_d\bigl(C_K+\|\omega\|_{\Dini}+C_{\WBP}+
        \norm{T1}_{\BMO}+\norm{T^*1}_{\BMO}\bigr),
\end{equation}
where $C_{\WBP}$ is the constant from the weak boundedness property. 
\end{theorem}

The main point of the proof is a decomposition
\begin{equation}\label{eq:intro-decomposition}
 T=S+E.
\end{equation}
The error $E$ is bounded directly on $L^2$, with a bound depending
only on the kernel. The smooth part $S$ is then treated by the
classical $T1$ theorem. The following result is the
precise form of this idea.

\begin{theorem}\label{thm:regularization}
Let $K$ satisfy \eqref{eq:size} and \eqref{eq:regularity}, with
$\omega$ satisfying the Dini condition. Then there exist an operator $E$ bounded on $L^2(\R^d)$
and a kernel $H\in C^\infty(\R^{d}\times \R^d\setminus\{(x, x): x\in \R^d\})$ such that
$K-H$ is the kernel of $E$,
\begin{equation}\label{eq:error-bound}
 \norm E_{L^2\to L^2}\le C_d(C_K+\|\omega\|_{\Dini}),
\end{equation}
and, for every pair of multiindices $\alpha,\beta\in \N^d$,
\begin{equation}\label{eq:smooth-bound}
 |\partial_x^\alpha\partial_y^\beta H(x,y)|
 \le C_{d,\alpha,\beta}(C_K+\|\omega\|_{\Dini})
       |x-y|^{-d-|\alpha|-|\beta|}.
\end{equation}
\end{theorem}

Let $T$ be as in Theorem~\ref{thm:main}, and define $S:=T-E$. It is easy to check that $H$ is the kernel of $S$. Weak boundedness property 
passes to $S$ immediately from \eqref{eq:error-bound}.
Moreover, by \eqref{eq:smooth-bound} the kernel of $E$, namely $K-H$, satisfies the Dini condition. Standard arguments therefore give
$E1,E^*1\in\BMO$, with norms at most $C_d(C_K+\|\omega\|_{\Dini})$.
Consequently $S$ inherits all the $T1$ assumptions of $T$.
Notice that here the kernel information on $E$ is also important, \eqref{eq:error-bound} alone would not imply the
required $\BMO$ statements.

We should mention that the idea of reducing a singular integral to a
simpler operator modulo a bounded error appears already in the literature.  For example, the
BCR-based construction of Auscher and Yang \cite{AY} produces
a perfect dyadic operator for H\"older kernels under an additional
local boundedness assumption on the full kernel.

Now we briefly explain how the error operator $E$ is found. First of all, we localize $K$ to
$|x-y|\simeq2^{-j}$, this means that we write 
\[
K(x,y)=\sum_j K_j(x,y),\qquad \supp K_j \subseteq \big\{(x,y): |x-y|\simeq2^{-j}\big\}.
\]
Denote by $T_j$ the integral operator associated with $K_j$. Let $P_j$ be a Fourier multiplier with symbol
$\rho (2^{-j} \xi)$, where $\rho\in C_c^\infty(\R^d)$ is a radial function such that $0\le\rho\le1$,
$\rho=1$ on $\{|\xi|\le1\}$, and $\rho=0$ on
$\{|\xi|\ge2\}$. Formally, we write 
\[
T_j= (I-P_j)T_j(I-P_j)+ P_j T_j(I-P_j)+ (I-P_j)T_jP_j + P_j T_j P_j.
\]
We shall prove that 
\[
\sum_j  (I-P_j)T_j(I-P_j)
\]
converges strongly in $L^2$ (Lemma \ref{lem:keykey}). 
The part
$
\sum_j (I-P_j)T_jP_j
$
needs an additional correction. We will introduce a smooth paraproduct $\Pi_{b_j}$ and show that 
\[
\sum_j  \big( (I-P_j)T_jP_j- \Pi_{b_j}\big)
\]
converges strongly in $L^2$ (Lemma \ref{lem:keykey}). Since 
\[
( (I-P_j)T_jP_j )^*= P_j T_j^* (1-P_j)
\]
this can be handled in a similar way. These three terms will form the error operator $E$. 

To get a heuristic feeling that why we are able to achieve the Dini regularity we also mention the following:  
\[
\sum_{q= 1}^\infty 2^q \omega (2^{-2^{q}})\lesssim \sum_{q= 1}^\infty \int_{2^{q-1}}^{2^q}\omega (2^{-t})
\ud t=\int_1^\infty \omega (2^{-t})
\ud t \simeq \int_0^1 \omega(t) \frac{\ud t}{t}.
\]
The sum in the left hand side of the above display is what we shall use, see the proof of  Theorem~\ref{thm:regularization}.

This paper is organized as follows. Section \ref{sec:preliminaries} 
is devoted to providing the details of our decomposition of the kernel, and some auxiliary results. Section \ref{sec:proof} is focused on the proof of the main results. 

\vspace{0.3cm}

\noindent \textbf{Acknowledgements.}
This work is supported by the National Natural Science Foundation of China through
project numbers 12671125 and 12222114. 

\vspace{0.3cm}

\noindent \textit{AI use statement.} 
The proof strategy is produced with the help of GPT-6-Astra ultra. The paper is written by the author, who will take the full responsibility.

\section{Preliminaries}\label{sec:preliminaries}
This section is devoted to recording basic notations and preliminary estimates, which will be needed in subsequent sections. 
\subsection{Classical singular integral theory}
As usual we use $\mathscr S(\R^d)$ to stand for the Schwartz class and 
denote its dual by
$\mathscr S'(\R^d)$. A singular integral operator is initially a continuous linear map
\[
 T:\mathscr S(\R^d)\longrightarrow\mathscr S'(\R^d).
\]
Its dual is defined by $\langle T^*f, g\rangle= \langle f, T g\rangle$. We say that $K: \R^d \times \R^d \setminus \{(x,x): x\in \R^d\}\to \C$ is the kernel of $T$ if we have the kernel representation 
\begin{equation}\label{pre:offdiagonal}
 \pair{Tf}{g}
 =\iint_{\R^d\times\R^d}K(x,y)f(y)g(x)\,dy\,dx
\end{equation}
whenever $f,g\in C_c^\infty(\R^d)$ have disjoint supports. Throughout the paper we assume $K$ satisfies 
the size condition \eqref{eq:size} and the regularity condition \eqref{eq:regularity} with an increasing
subadditive function $\omega:[0,1]\to[0,\infty)$ with $\omega(0)=0$ and 
\[
\|\omega\|_{\Dini}:=\int_0^1 \omega(t)\frac{\ud t}{t}<\infty. 
\]
In particular, subadditivity and monotonicity give
\begin{equation}\label{eq:sm}
\omega(st)\le(1+s)\omega(t)
\end{equation}
whenever $t, st \in [0,1]$, this fact will be frequently used later. 

Now let us state the $T1$ assumptions. 
Let $Q$ be a cube in $\R^d$ with sides parallel to the coordinate
axes, we call $\phi\in C_c^\infty(\R^d)$ an adapted normalized bump if
\[
 \supp \phi\subset Q,
 \qquad
 \ell(Q)^{|\alpha|}\norm{\partial^\alpha\phi}_\infty\le1
 \quad (|\alpha|\le d+2).
\]
We say that $T$ has the weak boundedness property if
\begin{equation}\label{pre:wbp}
 |\pair{T\phi}{\psi}|\le C_{\WBP} |Q|
\end{equation}
for every cube $Q$ and every pair of such bumps $\phi,\psi$.
This is the usual finite-order bump formulation. 
We also assume $T1, T^*1\in \BMO$, recall that 
\[
 \norm{b}_{\BMO}
 :=\sup_Q\frac1{|Q|}\int_Q|b-b_Q|,
 \qquad b_Q:=\frac1{|Q|}\int_Qb.
\]

Recall that in the introduction we have mentioned that the main goal is to prove the decomposition 
$T=S+E$, where $E$ is directly bounded on $L^2$ and we shall apply the classical $T1$ theorem to the smooth operator $S$. But the $T1$ assumptions are assumed on $T$, so we have to prove that $S$ can inherit the $T1$ assumptions from $T$. 

The weak boundedness property of $S$ is trivial. For the others we record the following 

\begin{lemma}\label{pre:bounded-bmo}
Let $E$ be a singular integral operator whose kernel satisfies \eqref{eq:size} and \eqref{eq:regularity} with a Dini modulus $\theta$. Suppose that $E$ is bounded on $L^2$, then
\[
 \norm{E1}_{\BMO}+\norm{E^* 1}_{\BMO}
 \le C\bigl(\norm{E}_{L^2\to L^2}+\norm{\theta}_{\Dini}\bigr).
\]
\end{lemma}
The proof is standard, since Dini condition implies the classical H\"ormander condition, in which case we already know that singular integrals maps $L^\infty$ to $\BMO$. Next we record the classical 
David--Journ\'e $T1$ theorem.

\begin{theorem}\label{pre:classical-t1}
Let $S:\mathscr S(\R^d)\to\mathscr S'(\R^d)$ be continuous and the associated
 kernel $H$  satisfies
\[
 |H(x,y)|\le C_H|x-y|^{-d},\qquad
 |\nabla_xH(x,y)|+|\nabla_yH(x,y)|
       \le C_H|x-y|^{-d-1}.
\]
If $S$ satisfies \eqref{pre:wbp}  and 
$S1,S^*1\in\BMO$, then $S$ has a unique bounded extension to
$L^2(\R^d)$, and
\[
 \norm{S}_{L^2\to L^2}
 \le C_d\bigl(C_H+C_\WBP+
        \norm{S1}_{\BMO}+\norm{S^*1}_{\BMO}\bigr).
\]
\end{theorem}

\subsection{Annular decomposition of the kernel}
Now we move to the details of the annular decomposition of the kernel. 
Choose a smooth radial function $\chi_0$ such that $\chi_0=1$
on $\{|u|\le1\}$ and $\chi_0=0$ on $\{|u|\ge2\}$, and set
$\eta(u)=\chi_0(u)-\chi_0(2u)$. Then $\eta$ is supported in
$\{1/2\le|u|\le2\}$ and telescoping gives
\[
 \sum_{j\in\mathbb Z}\eta(2^ju)=1,\qquad u\ne0.
\]
Then  define
\begin{equation}\label{pre:annuli}
 K_j(x,y)=\eta(2^j (x-y) )K(x,y),\qquad
 T_jf(x)=\int_{\R^d}K_j(x,y)f(y)\,dy.
\end{equation}
Each $T_j$ is a genuine integral
operator, independently of whether $T$ is bounded. The size
estimate \eqref{eq:size} and the annular support give
\[
|T_j f| \lesssim C_K Mf,
\]
where $M$ stands for the Hardy-Littlewood maximal operator. Hence we easily have $\norm{T_j}_{L^2\to L^2}\le C_K$.
For every $x\ne y$,
\begin{equation*} 
 K(x,y)=\sum_{j\in\mathbb Z}K_j(x,y),
\end{equation*}
and this sum has only finitely many nonzero terms whenever $x\ne y$. It follows that 
if $f,g\in C_c^\infty(\R^d)$ have disjoint supports, then 
\begin{equation}\label{pre:separated-sum}
 \pair{Tf}{g}=\sum_{j\in\mathbb Z}\pair{T_jf}{g}
 \qquad
 (\supp f\cap\supp g=\varnothing),
\end{equation}
where the sum in this formula is finite. 
Next we conduct the further decomposition of $T_j$ 
mentioned in the introduction in details. 

Let $P_j$ be a Fourier multiplier with symbol
$\rho (2^{-j} \xi)$, where $\rho\in C_c^\infty(\R^d)$ is a radial function such that $0\le\rho\le1$,
$\rho=1$ on $\{|\xi|\le1\}$, and $\rho=0$ on
$\{|\xi|\ge2\}$. 
Define  $Q_j=P_j-P_{j-1}$. We also define the related non-smooth counterpart. Namely, let 
$L_j$ be the ball Fourier multiplier with symbol $1_{\{|\xi|\le 2^j\}}$ and define $H_j=I-L_j$, here $I$ stands for the identity operator. 
By definition, 
\begin{equation*} 
 \supp\widehat{Q_j f}
 \subset\{\xi:2^{j-1}\le|\xi|\le2^{j+1}\}.
\end{equation*}
Then  Plancherel's theorem
and bounded overlap give
\begin{equation}\label{pre:q-square}
 \sum_{j\in\mathbb Z}\norm{Q_j f}_2^2\le C\norm f_2^2.
\end{equation}
 With all these at hand, we decompose 
 \[
T_j= (I-P_j)T_j(I-P_j)+ (I-P_j)T_jP_j  + P_j T_j(I-P_j)+ P_j T_j P_j.
\]
Since $T_j, P_j$ are well-defined on $L^2$ functions, all terms on the right hand side are well-defined. 
The following result will be frequently used in order to deal with these terms. 
\begin{proposition}\label{prop:keytool}
Let $F_j$ be a sequence of bounded operators on $L^2(\R^d)$ with $\|F_j\|_{L^2\to L^2}\le C_0$. Suppose that for every integer $N\ge 10$, 
\[
\| (I-B_{j, N}) F_j \|_{L^2 \to L^2}+ \| F_j  (I-B_{j,N})\|_{L^2 \to L^2}\le a_N,
\]
where $B_{j,N}= L_{j+N}-L_{j-N}$ and $a_N\ge 0$ is decreasing and summable over $N\ge 10$. Then 
\[
\Big\| \sum_j F_j \Big\|_{L^2 \to L^2}\le CC_0+ C\sum_{N\ge 6}a_N.
\]
\end{proposition}
\begin{proof}
Since  
\begin{align*}
\| F_j - B_{j,N} F_j B_{j,N}\|_{L^2 \to L^2}&\le \| F_j - B_{j,N} F_j\|_{L^2 \to L^2}+\| B_{j,N} F_j- B_{j,N} F_j B_{j,N}\|_{L^2 \to L^2}\\
&\le \| (I-B_{j,N}) F_j \|_{L^2 \to L^2}+ \| F_j (I- B_{j,N}) \|_{L^2 \to L^2}\le a_N,
\end{align*}
and $a_N \to 0$ as $N\to \infty$, we can write 
\begin{align*}
 F_j =  B_{j,10}F_j B_{j,10}+ \sum_{q=1}^\infty \big(B_{j,10\cdot 2^q}F_j B_{j,10\cdot 2^q} 
 - B_{j,10\cdot 2^{q-1}}F_j B_{j,10\cdot 2^{q-1}} \big).
 \end{align*}
 Then for any finite $\mathcal J \subset \Z$ we have 
 \begin{align*}
 \Big\|\sum_{j\in \mathcal J}  B_{j,10}F_j B_{j,10}f\Big\|_{L^2}
 &\le \sum_{r=0}^{19}   \Big\|\sum_{\substack{j\in \mathcal J \\ j\equiv r({\rm mod} \,20)}}  B_{j,10}F_j B_{j,10}f\Big\|_{L^2}.
 \end{align*}
By plancherel, and set the Fourier support of $B_{j,10} $ by $U_{j, 10}$, where we define
\[
U_{j, N}:=\{\xi: 2^{j-N} \le |\xi|\le 2^{j+N}\},
\] we have
 \begin{align*}
\Big\|\sum_{\substack{j\in \mathcal J \\ j\equiv r({\rm mod} \,20)}}  B_{j,10}F_j B_{j,10}f\Big\|_{L^2}
&=\Big( \sum_{\substack{j\in \mathcal J \\ j\equiv r({\rm mod} \,20)}}\big\|   B_{j,10}F_j B_{j,10}f\big\|_{L^2}^2 \Big)^{1/2}\\
&\le  C_0 \Big( \sum_{\substack{j\in \mathcal J \\ j\equiv r({\rm mod} \,20)}}\big\|   B_{j,10}f\big\|_{L^2}^2 \Big)^{1/2}\\
&\le C_0 \|\widehat f \, 1_{\cup_{j\in \mathcal J}U_{j, 10}}\|_{L^2}\le C_0 \|f\|_{L^2}.
 \end{align*}
This  shows that 
 \[
 \Big\|\sum_{j\in \Z}  B_{j,10}F_j B_{j,10}f\Big\|_{L^2}\lesssim C_0\|f\|_{L^2}.
 \]
 
 Next we control the norm of
\begin{equation}\label{eq:normconv}
\sum_{j\in \mathcal J} \sum_{q=1}^\infty \big(B_{j,10\cdot 2^q}F_j B_{j,10\cdot 2^q} 
 - B_{j,10\cdot 2^{q-1}}F_j B_{j,10\cdot 2^{q-1}} \big),
\end{equation}
where again $\mathcal J$ is a finite set. 
Then we can change the order of the summation and we can rewrite the above as 
\begin{align*}
&\sum_{q=1}^\infty\sum_{j\in \mathcal J}  \big(B_{j,10\cdot 2^q}F_j B_{j,10\cdot 2^q} 
 - B_{j,10\cdot 2^{q-1}}F_j B_{j,10\cdot 2^{q-1}} \big)\\
 &= \sum_{q=1}^\infty \sum_{r=0}^{20\cdot 2^q-1}\sum_{\substack{j\in \mathcal J\\ j\equiv r({\rm mod} \,20\cdot 2^q)}  } \big(B_{j,10\cdot 2^q}F_j B_{j,10\cdot 2^q} 
 - B_{j,10\cdot 2^{q-1}}F_j B_{j,10\cdot 2^{q-1}} \big).
\end{align*}
Similar as before we obtain 
\begin{align*}
&\Big\|\sum_{q=1}^\infty\sum_{j\in \mathcal J}  \big(B_{j,10\cdot 2^q}F_j B_{j,10\cdot 2^q} 
 - B_{j,10\cdot 2^{q-1}}F_j B_{j,10\cdot 2^{q-1}} \big)f\Big\|_{L^2}\\
 &\le  \sum_{q=1}^\infty \sum_{r=0}^{20\cdot 2^q-1}
\Big( \sum_{\substack{j\in \mathcal J\\ j\equiv r({\rm mod} \,20\cdot 2^q)}  } 
 \big\|\big(B_{j,10\cdot 2^q}F_j B_{j,10\cdot 2^q} 
 - B_{j,10\cdot 2^{q-1}}F_j B_{j,10\cdot 2^{q-1}} \big)f\big\|_{L^2}^2\Big)^{\frac 12}.
\end{align*}
At this point let us observe that 
\begin{align*}
&\big(B_{j,10\cdot 2^q}F_j B_{j,10\cdot 2^q} 
 - B_{j,10\cdot 2^{q-1}}F_j B_{j,10\cdot 2^{q-1}}  \big)f\\
 &\hspace{2cm} =\big(B_{j,10\cdot 2^q}F_j B_{j,10\cdot 2^q} 
 - B_{j,10\cdot 2^{q-1}}F_j B_{j,10\cdot 2^{q-1}} \big)(B_{j,10\cdot 2^q }f)
\end{align*}
and that 
\begin{align*}
&\|B_{j,10\cdot 2^q}F_j B_{j,10\cdot 2^q} 
 - B_{j,10\cdot 2^{q-1}}F_j B_{j,10\cdot 2^{q-1}} \|_{L^2\to L^2}\\
 &\le \| B_{j,10\cdot 2^q}F_j B_{j,10\cdot 2^q} 
 - F_j \|_{L^2\to L^2}+ \| F_j
 - B_{j,10\cdot 2^{q-1}}F_j B_{j,10\cdot 2^{q-1}}   \|_{L^2\to L^2}\\
 &\le a_{10\cdot 2^q} +a_{10\cdot 2^{q-1}}\le 2a_{10\cdot 2^{q-1}}.
\end{align*}
Combining this we get
\begin{align*}
&\Big\|\sum_{q=1}^\infty\sum_{j\in \mathcal J}  \big(B_{j,10\cdot 2^q}F_j B_{j,10\cdot 2^q} 
 - B_{j,10\cdot 2^{q-1}}F_j B_{j,10\cdot 2^{q-1}} \big)f\Big\|_{L^2 \to L^2}\\
 &\le
C \sum_{q=1}^\infty  2^q a_{10\cdot 2^{q-1}}\|\wh f \,1_{\cup_{j\in \mathcal J} U_{j, 10\cdot 2^q}}\|_{L^2}.
\end{align*}
This shows that \eqref{eq:normconv} converges and the limit is bounded on $L^2$ with norm 
\[
 \sum_{q=1}^\infty  2^q a_{10\cdot 2^{q-1}}\le \frac 25\sum_{q=1}^\infty  \sum_{s=10\cdot 2^{q-2}+1}
 ^{10\cdot 2^{q-1}}a_s=\frac 25\sum_{q\ge 6}^\infty a_q.
 \]
This completes the proof. 
\end{proof}
Proposition \ref{prop:keytool} is crucial for us. We will apply it to 
$F_j=(I-P_j)T_j(I-P_j) $ and $F_j=(I-P_j)T_jP_j -\Pi_{b_j}$, where $b_j= (I-P_j)T_j1$ and 
\[
\Pi_{b_j}f:= \sum_{n\in\mathbb Z}(Q_nb_j)P_{n-6}f.
\]
We provide some observations in order. Firstly, since $T_j1$ is bounded, $b_j$ is bounded as well. Secondly, the definition of $Q_n$ and $b_j$ implies that $Q_n b_j=0$ whenever $n\le j-1$. Lastly, whenever $f$ has Fourier support at $\{\xi: |\xi|\le 2^{j-N}\}$ for $N\ge 6$, we have 
\begin{equation}\label{eq:ob}
\Pi_{b_j}f= \sum_{n\ge j}(Q_nb_j)P_{n-6}f=\sum_{n\ge j}(Q_nb_j) f= b_j f.
\end{equation}
\eqref{eq:ob} explains why we do such a correction, since in this way we get some cancellation in the kernel of $(I-P_j)T_jP_j- \Pi_{b_j}$ when hitting on a function with suitable Fourier support. Note that  $P_j T_j(I-P_j)$ can be viewed as the dual of $(1-P_j)T_j^* P_j$, so it can be handled in a similar way. And 
$P_j T_j P_j$ is a component of the smooth part, we will discuss this in the next section.

We end up this section with the following (well-known) lemma.

\begin{lemma}\label{pre:paraproduct}
Let $b\in L^\infty$ and $\Pi_{b}$ be defined as \eqref{eq:ob}. Then we have 
\begin{equation}\label{pre:para-bound}
 \norm{\Pi_{b}}_{L^2\to L^2}\le C \|b\|_{L^\infty}.
\end{equation}
Moreover, the associated kernel 
$
\pi(x,y)$ 
satisfies that for all multi-indices $\alpha,\beta$,
\begin{equation}\label{pre:para-kernel-derivative}
 |\partial_x^\alpha\partial_y^\beta\pi (x,y)|
 \le C_{\alpha,\beta} \|b\|_{L^\infty}
       |x-y|^{-d-|\alpha|-|\beta|},\qquad x\ne y.
\end{equation}
\end{lemma}

\begin{proof}
The estimate \eqref{pre:para-bound} is well-known, for completeness we sketch the proof here. 
First of all, note that the Fourier support of the summand $(Q_n b ) P_{n-6} f $ is 
\[
\big\{\xi: (2^{-1}-2^{-5})2^n\le |\xi|\le (2+2^{-5}) 2^n\big\},
\]
then we may split the paraproduct as 
\[
\Pi_{b }f= \sum_{r=0}^4 \sum_{n\equiv r\,(\rm{mod}5)} (Q_n b ) P_{n-6}f=: \sum_{r=0}^4 \Pi_{b }^r f.
\]
By triangle's inequality we only need to bound $\Pi_{b }^r f$, and now the summands in $\Pi_{b }^r f$ have disjoint Fourier support,
using Plancherel's formula we only need to estimate the diagonal 
\[
\sum_{n\equiv r\,(\rm{mod}5)} \int_{\R^d }|Q_n b |^2 |P_{n-6}f|^2
= \sum_{n\equiv r\,(\rm{mod}5)} \sum_{Q\in \calD_n}\int_Q |Q_n b |^2 |P_{n-6}f|^2,
\]
where $\calD_n$ stands for the set of standard dyadic cubes with side-length $2^{-n}$. On one hand, for $x\in Q\in \calD_n$,
\begin{align*}
|P_{n-6}f(x)|&\lesssim \int_{\R^d} \frac {2^{(n-6)d}}{(1+2^{n-6}|x-y|)^{100d}}|f(y)|\ud y\\
&\lesssim \inf_{z\in Q} \int_{\R^d} \frac {2^{(n-6)d}}{(1+2^{n-6}|z-y|)^{100d}}|f(y)|\ud y\lesssim \inf_{z\in Q} Mf(z).  
\end{align*}
On the other hand, it is not difficult to establish the Carleson estimate
\begin{equation}\label{eq:para-carleson}
\sum_{n\equiv r\,(\rm{mod}5)}\sum_{\substack{Q\in \calD_n \\ Q\subset Q_0 }} \int_Q |Q_n b |^2\le \sum_{n\equiv r\,(\rm{mod}5)}  \int_{Q_0} |Q_n b |^2\lesssim \|b\|_{L^\infty}^2 |Q_0|,
\end{equation}
where $r\in \{0,1,\ldots, 4\}$ is fixed. 
In fact, split $b =b \one_{4Q_0}+b \one_{(4Q_0)^c}$, where $4Q_0$ is the
concentric dilation of $Q_0$. For the first part, using
\eqref{pre:q-square} we have 
\[
\sum_{n\equiv r\,(\rm{mod}5)}\sum_{\substack{Q\in \calD_n \\ Q\subset Q_0 }} \int_Q |Q_n (b 1_{4Q_0})|^2\le  \sum_n\norm{Q_n(b \one_{4Q_0})}_2^2
 \le C\norm{b \one_{4Q_0}}_2^2
 \le C \|b\|_{L^\infty}^2|Q_0|.
\]
For the second part, using the fast decay of the kernel (of $Q_n$, say $q_n$) we have 
\[
 |Q_n(b \one_{(4Q_0)^c})(x)|
 \le C  \|b\|_{L^\infty} (2^n\ell(Q_0))^{-99d},
 \qquad x\in Q_0,\quad n\ge n_0,
\]where $n_0$ is the generation of $Q_0$, i.e. we assume $Q_0\in \calD_{n_0}$. 
Indeed, integrate the estimate
$|q_n(z)|\le C 2^{nd}(1+2^n|z|)^{-100d}$
over $|z|\ge c\ell(Q)$. Squaring and summing the resulting
geometric series proves \eqref{eq:para-carleson}.

Now with \eqref{eq:para-carleson} at hand, we have 
\begin{align*}
\sum_{n\equiv r\,(\rm{mod}5)} \sum_{Q\in \calD_n}\int_Q |Q_n b |^2 |P_{n-6}f|^2
&\lesssim \sum_{n\equiv r\,(\rm{mod}5)} \sum_{Q\in \calD_n}\int_Q |Q_n b |^2 (\inf_Q Mf)^2 \\
&= 2\int_0^\infty  \sum_{n\equiv r\,(\rm{mod}5)} \sum_{\substack{Q\in \calD_n\\ \inf_Q Mf >\lambda }}\int_Q |Q_n b |^2 \lambda\ud \lambda\\
&\le 2\int_0^\infty  \sum_{n\equiv r\,(\rm{mod}5)} \sum_{\substack{Q\in \calD_n\\ Q\subset \{ Mf >\lambda\} }}\int_Q |Q_n b |^2 \lambda\ud \lambda\\
&\lesssim \|b\|_{L^\infty}^2 \int_0^\infty  \lambda |\{Mf>\lambda\}|\ud \lambda\\
&\simeq  \|b\|_{L^\infty}^2\|Mf\|_{L^2}^2 \lesssim \|b\|_{L^\infty}^2 \|f\|_{L^2}^2.
\end{align*}

It remains to prove \eqref{pre:para-kernel-derivative}. 
Notice that every
differentiated summand of $\pi (x,y)$ can be bounded by 
\begin{align*}
 C_{\alpha,\beta,N}\| b  \|_{L^\infty}
 2^{n(d+|\alpha|+|\beta|)}(1+2^n|x-y|)^{-N}.
\end{align*}
Then simply write 
\[
\sum_n = \sum_{n: 2^n \ge |x-y|^{-1}}+  \sum_{n: 2^n < |x-y|^{-1}}
\]
gives the desired estimate. 
\end{proof}

\section{Proof of the main results}\label{sec:proof}

We first prove Theorem~\ref{thm:regularization}. We begin with the following kernel estimates
\begin{lemma} 
Denote the kernels of $T_j P_j$ and $P_j T_j P_j$ by $K_j^{(1)}$ and $K_j^{(2)}$, respectively.
Then there holds that 
\begin{equation}\label{eq:K1}
|K_j^{(1)}(x,y)|\le C C_K 2^{jd}(1+2^j |x-y|)^{-100d}
\end{equation}
and for all multi-indices $\alpha, \beta\in \N^d$,
\begin{equation}\label{eq:K2}
|\partial_x^\alpha \partial_y^\beta K_j^{(2)}(x,y)|\le C_{\alpha, \beta} C_K 2^{j(d+|\alpha|+|\beta|)}(1+2^j |x-y|)^{-100d}.
\end{equation}
\end{lemma}
\begin{proof}
We first prove \eqref{eq:K1}. By definition, we can write 
\[
K_j^{(1)}(x,y):= \int K_j(x,v) p_j(v-y) \ud v,
\]
where $p_j$ is defined via $\wh {p_j}(\xi)= \rho(2^{-j}\xi)$. Then 
since $|x-v|\sim 2^{-j}$ we have 
\[
|p_j(v,y)|\lesssim 2^{jd} (1+2^j|v-y|)^{-100d}\sim 2^{jd} (1+2^j|x-y|)^{-100d}. 
\]
Hence \eqref{eq:K1} follows. Next we prove \eqref{eq:K2}. Again the kernel of  $P_j T_j P_j$ can be written as 
\[
  K_j^{(2)}(x,y)=\iint p_j(x-u)K_j(u,v)p_j(v-y)\ud u \ud v.
\]
Then 
\[
\partial_x^\alpha \partial_y^\beta K_j^{(2)} (x, y) = 2^{j( |\alpha|+|\beta|)}\iint 2^{jd}(\partial^\alpha p)(2^j(x-u))K_j(u,v)2^{jd}(\partial^\beta p)(2^j(v-y))\ud u \ud v.
\]
As before we have 
\[
|(\partial^\alpha p)(2^j(x-u))|\lesssim (1+ 2^j |x-u|)^{-100d}
\]
and
\[
|(\partial^\beta p)(2^j(v-y))|\lesssim (1+ 2^j |v-y|)^{-100d}.
\]
Since $|u-v|\sim 2^{-j}$, we have 
\[
1+ 2^j |x-u|\sim 1+ 2^j |x-v|.
\]
Then split the integration over $v$ into $|x-v|\ge |x-y|/2$ and $|x-v|<|x-y|/2$, in the latter case we must have 
$|y-v|\ge |x-y|/2$. Hence we conclude that 
\[
|\partial_x^\alpha \partial_y^\beta K_j^{(2)}(x, y) |\lesssim C_K 2^{j( |\alpha|+|\beta|)} 2^{jd}(1+ 2^j |x-y|)^{-100d}.
\]
\end{proof} 
Next we prove the following, which is crucial to the main results. 
\begin{lemma}\label{lem:keykey}
Let 
\[
F_j\in \big\{ (I-P_j)T_j(I-P_j),  (I-P_j)T_jP_j- \Pi_{b_j}\big\}.
\] Then there holds that  
\begin{equation*} 
 \Big\| \sum_j F_j  \Big\|_{L^2 \to L^2 }\le C(\|\omega\|_{\Dini}+C_K).
 \end{equation*}
\end{lemma}
\begin{proof}
With Proposition \ref{prop:keytool} at hand, it suffices to prove that 
for every integer $N\ge 10$, 
\[
\| (I-B_{j, N}) F_j \|_{L^2 \to L^2}+ \| F_j  (I-B_{j,N})\|_{L^2 \to L^2}\le C\nu(C2^{-N}), 
 \]
 where 
$
 \nu(t)= \omega(\min\{ct,1\})+C_Kt 
$, and $c, C$ are some absolute constants. 
 
We first consider the case $F_j= (I-P_j)T_j(I-P_j)$. Note that in this case we have 
\[
(I-B_{j, N}) F_j= H_{j+N}F_j = (I-P_j) H_{j+N} T_j (I-P_j)
\]
and 
\[
 F_j(I-B_{j, N})= F_jH_{j+N}= (I-P_j)T_j H_{j+N} (I-P_j).
\]
Therefore, it suffices to prove \[
\| H_{j+N}  T_j \|_{L^2} + \| T_j H_{j+N} \|_{L^2} \le C\nu(C2^{-N}).
\]
We only focus on $H_{j+N}  T_j $ since the other is similar (note that $(T_j H_{j+N})^*= H_{j+N} T_j^*$). 
Now a crucial trick is, by Plancherel we have
\[
\| H_{j+N}  T_j \|_{L^2} \le \| (1- P_{j+N-1})T_j \|_{L^2}.
\]
This is important since the kernel of $(1- P_{j+N-1})T_j $ is 
\[
\int (K_j (x,y) - K_j (x-h, y) ) 2^{(j+N-1) d} (\mathcal F^{-1} \rho) ( 2^{j+N-1}h ) \ud h, 
\]
where $\mathcal F^{-1} $ stands for the inverse Fourier transform, the smooth cut-off has made it easier to handle. 
By Schur's test it suffices to control 
\begin{align*}
&\sup_x \int\Big|\int (K_j (x,y) - K_j (x-h, y) ) 2^{(j+N-1) d} (\mathcal F^{-1} \rho) ( 2^{j+N-1}h ) \ud h\Big|\ud y\\
&\hspace{1cm}\le \int \sup_x\int  |K_j (x,y) - K_j (x-h, y) |\ud y 2^{(j+N-1) d} |(\mathcal F^{-1} \rho) ( 2^{j+N-1}h )| \ud h
\end{align*}
and 
\begin{align*}
&\sup_y \int\Big|\int (K_j (x,y) - K_j (x-h, y) ) 2^{(j+N-1) d} (\mathcal F^{-1} \rho) ( 2^{j+N-1}h ) \ud h\Big|\ud x\\
&\hspace{1cm}\le \int \sup_y\int  |K_j (x,y) - K_j (x-h, y) |\ud x 2^{(j+N-1) d} |(\mathcal F^{-1} \rho) ( 2^{j+N-1}h )| \ud h.
\end{align*}
First of all, write 
\begin{align*}
|K_j(x,y)- K_j(x-h, y)|&= |\eta(2^j (x-y) )K(x,y) - \eta(2^j (x-y-h) )K(x-h,y)|\\
&\le  |\eta(2^j (x-y) ) - \eta(2^j (x-y-h) ) |\cdot |K(x,y)|\\
&\qquad+ |\eta(2^j (x-y-h) ) |\cdot |K(x,y)-K(x-h, y)|.
\end{align*}
Whenever $|h|\le 2^{-j-5}$, if $\eta(2^j (x-y) ) - \eta(2^j (x-y-h) )$ is non-zero then we must have 
$|x-y|\sim 2^{-j}$, then 
\[
 |\eta(2^j (x-y) ) - \eta(2^j (x-y-h) ) |\cdot |K(x,y)|\lesssim C_K2^j |h| 2^{jd}.
\]
On the other hand, $ \eta(2^j (x-y-h) )$ is non-zero also implies $|x-y|\sim 2^{-j}$,
 by 
the regularity of the kernel \eqref{eq:regularity}, we get 
\[
 |\eta(2^j (x-y-h) ) |\cdot |K(x,y)-K(x-h, y)|\lesssim   \omega (c2^j  |h|)2^{jd}.
\]
As a consequence we have 
\[
\int |K_j(x, y)- K_j (x-h, y)|\ud y + \int |K_j(x, y)- K_j (x-h, y)|\ud x\lesssim \nu (2^j |h|). 
\]
Now consider $|h|> 2^{-j-5}$, then by the size estimate \eqref{eq:size} we also have
\[
\int |K_j(x, y)- K_j (x-h, y)|\ud y + \int |K_j(x, y)- K_j (x-h, y)|\ud x\lesssim C_K\lesssim \nu (2^j |h|).
\]
Then we see that both objects in Schur test are bounded by 
\begin{align*}
&\int 2^{(j+N-1) d} \nu (2^j |h|)|(\mathcal F^{-1} \rho) ( 2^{j+N-1}h )| \ud h\\
&=  \int \nu(2^{1-N} |h|)|(\mathcal F^{-1} \rho) (  h )| \ud h\le \nu(2^{1-N}) \int (|h|+1)|(\mathcal F^{-1} \rho) (  h )| \ud h,
\end{align*}
where we have used \eqref{eq:sm}. This finishes the case $F_j=(I-P_j)T_j (I-P_j)$. 

Next let us turn to the case $F_j=(I-P_j)T_jP_j- \Pi_{b_j}$. Again, observe the following 
\[
(I- B_{j,N}) F_j = (I-P_j) H_{j+N} T_jP_j- H_{j+N} \Pi_{b_j}
\]
and 
\[
F_j(I- B_{j,N})= \big((I-P_j)T_j L_{j-N}- \Pi_{b_j} L_{j-N}\big) - \Pi_{b_j} H_{j+N}.
\]
It is clear that 
\[
\|(I-P_j) H_{j+N} T_jP_j\|_{L^2\to L^2}\le \|  H_{j+N} T_j\|_{L^2\to L^2}\le C\nu(C2^{-N}).
\]
On the other hand, by the Fourier support of $(Q_n b_j)P_{n-6}f $, we know that 
$H_{j+N}((Q_n b_j)P_{n-6}f)$ is nonzero unless $n\ge j+N-2$. Hence 
\[
H_{j+N} \Pi_{b_j}= H_{j+N} \Pi_{(1-P_{j+N-4})b_j}
\]
and it follows that 
\[
\|H_{j+N} \Pi_{b_j}\|_{L^2\to L^2 }\le \|\Pi_{(I-P_{j+N-4})b_j}\|_{L^2\to L^2 }\le \|(I-P_{j+N-4})b_j\|_{L^\infty}=  \|(I-P_{j+N-4})T_j1\|_{L^\infty}.
\]
Similar as above (the estimate of the kernel of $(I-P_{j+N-1})T_j$ we can conclude that 
\[
\|H_{j+N} \Pi_{b_j}\|_{L^2\to L^2 }\le C\nu(C2^{-N}).
\]
This finishes the estimate of $(I- B_{j,N}) F_j$.
Then note that in a similar way, $P_{n-6}H_{j+N}\ne 0$ implies that $n\ge j+N+5$. Hence 
\[
\Pi_{b_j}H_{j+N}= \Pi_{(1-P_{j+N+3})b_j}H_{j+N}
\]
and we still have 
\[
\| \Pi_{b_j}H_{j+N}\|_{L^2\to L^2 }\le C\nu(C2^{-N}).
\]
It remains to bound 
\[
\|(I-P_j)T_j L_{j-N}- \Pi_{b_j} L_{j-N}\|_{L^2\to L^2}=\|(I-P_j)T_jP_j L_{j-N}- \Pi_{b_j} L_{j-N}\|_{L^2\to L^2} .
\]
Recall that 
$
\Pi_{b_j} L_{j-N}f= b_jL_{j-N} f,
$ we reduce the estimate to the $L^2$ norm of 
\[
\int (K_j^{(1)}(x, y)-K_j^{(2)}(x, y))(L_{j-N}f(y)- L_{j-N}f(x))\ud y, 
\]
where recall that $K_j^{(1)}-K_j^{(2)}$ is the kernel of $(I-P_j)T_jP_j$, using \eqref{eq:K1} and 
\eqref{eq:K2} we have 
\[
|K_j^{(1)}(x, y)-K_j^{(2)}(x, y)|\le C C_K 2^{jd} (1+2^j|x-y|)^{-100d}.
\]
Hence by change of variable $y\to x-h$ and then using Minkowski's inequality we have that the object is controlled by 
\begin{align*}
 \int C C_K 2^{jd} (1+2^j|h|)^{-100d}\Big(\int    | L_{j-N}f(x-h)- L_{j-N}f(x)|^2 \ud x \Big)^{1/2}\ud h
\end{align*}
By Plancherel 
\[
\| L_{j-N}f(x-h)- L_{j-N}f(x)\|_{L^2}\le \sup_{|\xi|\le 2^{j-N}} |e^{-2\pi i h\cdot \xi}-1| \|f\|_{L^2}
\le 2^{j-N+1}\pi |h| \|f\|_{L^2}.
\]
Then we finally obtain 
\[
\|(I-P_j)T_jP_j L_{j-N}- \Pi_{b_j} L_{j-N}\|_{L^2\to L^2}\le C C_K 2^{-N} \|f\|_{L^2}. 
\]
This concludes the proof. 
\end{proof}

Now we are ready to prove Theorem~\ref{thm:regularization}. 
\begin{proof}[Proof of Theorem~\ref{thm:regularization}]
As we have indicated above, in our decomposition
\[
E= \sum_j (I-P_j) T_j (I-P_j) + \sum_j \big( (I-P_j) T_j P_j- \Pi_{b_j} \big)+ \sum_j \big( P_jT_j(I- P_j)- \Pi_{\wt b_j}^* \big),
\]
where $\wt b_j:= (I-P_j) T_j^* 1$. 
Lemma \ref{lem:keykey} already gives the boundedness of the first two terms on the right hand side, while the third term is actually similar since 
\[
P_jT_j(I- P_j)- \Pi_{\wt b_j}^*=\big( (I-P_j)T_j^*P_j- \Pi_{\wt b_j} \big)^*.
\] 
It remains to prove the kernel estimates. This will be proved in an indirect way. In fact, 
set $S=T-E$, 
when $f, g\in C_c^\infty$ with disjoint support, 
we have 
\begin{align*}
\langle Tf, g\rangle =\iint K(x,y) f(y) g(x) \ud y \ud x&=\Big\langle \sum_j T_j f,g\Big\rangle \\
&=\langle Sf, g\rangle+ \langle Ef, g\rangle.
\end{align*}
Then clearly we have 
\[
\langle Sf, g\rangle=\Big\langle \sum_j (P_j T_j P_j + \Pi_{b_j}+ \Pi_{\wt b_j}^* )f, g\Big\rangle.
\]
By \eqref{eq:K2}
we see that $\sum_j P_j T_j P_j$ has a smooth Calder\'on-Zygmund kernel. However, 
Lemma \ref{pre:paraproduct} only tells us that each individual term $ \Pi_{b_j}+ \Pi_{\wt b_j}^* $ has smooth Calder\'on-Zygmund kernel, but unfortunately we do not have 
the summability of $\sum_j \|b_j\|_{L^\infty}$. So we need some justifications here. In fact, 
since the kernel of $ \Pi_{b_j}$ is 
\[
\pi_j(x,y)= \sum_{n\ge j} (Q_n b_j)(x) p_{n-6}(x-y),
\]
where recall that $p_n$ is defined via $\wh {p_n}(\xi)= \rho(2^{-n}\xi)$. For any finite $\calJ$, we have 
\[
\pi_{\calJ}(x,y):=\sum_{j\in \calJ }\pi_j(x,y)= \sum_n \sum_{\substack{j\in \calJ \\ j\le n}}(Q_n b_j)(x) p_{n-6}(x-y).
\]
Crucially, in the proof of Lemma \ref{lem:keykey} we have proved that 
\[
\| T_j1 (\cdot-h)-T_j 1\|_{L^\infty}\le C \nu (2^j |h|).
\]
Then using the cancellation of the kernel of $Q_n$ and \eqref{eq:sm} we get that 
\[
\|Q_n b_j\|_{L^\infty}\le C \nu (2^{j-n}),
\]
where $n\ge j+10$. For $j\le n <j+10$ we use the elementary estimate 
\[
\|Q_n b_j\|_{L^\infty}\le C C_K. 
\]
Hence for any $n$ we have 
\[
\sum_{\substack{j\in \calJ \\ j\le n}}\| Q_n b_j\|_{L^\infty}\le  CC_K+ \sum_{\substack{j\in \calJ \\ j\le n}}  C\nu (2^{j-n}). 
\]
Since $Q_n b_j$ has Fourier support in an annulus of radius comparable to $2^n$, Bernstein's inequality implies
\[
\sum_{\substack{j\in \calJ \\ j\le n}}\| \partial^\alpha Q_n b_j\|_{L^\infty}\le 2^{n|\alpha|}\Big(C_\alpha C_K+ \sum_{\substack{j\in \calJ \\ j\le n}}  C_\alpha\nu (2^{j-n})\Big).
\]
So that we have 
\begin{align*}
|\partial^\alpha_x \partial_y \pi_{\calJ}(x,y)|&\lesssim \sum_n 2^{n(d+ |\alpha|+|\beta|)}
\Big(C_\alpha C_K+ \sum_{\substack{j\in \calJ \\ j\le n}}  C_\alpha\nu (2^{j-n})\Big)
(1+2^n|x-y|)^{-100d-|\alpha|-|\beta|}\\
&\le (C_K +\|\omega\|_{\Dini})|x-y|^{-d-|\alpha|-|\beta|}.
\end{align*}
Since we have actually proved the absolute convergence of the double sum, the same estimate holds for $\calJ=\Z$ as well. The kernel estimates for $\Pi_{\wt {b_j}}^*$ is similar. 
 This shows that $S$ has smooth Calder\'on-Zygmund kernel, and hence $E$ inherits the kernel estimates from $T$. 
\end{proof}

We now prove Theorem~\ref{thm:main}. 
\begin{proof}[Proof of Theorem~\ref{thm:main}]
The necessity is by-now standard and we omit the details. 

For sufficiency, with Theorem~\ref{thm:regularization} at hand, we see that 
$S=T-E$ inherits the weak boundedness property. Moreover,  standard Calder\'on-Zygmund theory gives that $E1, E^*1\in \BMO$, hence
$S1, S^*1\in \BMO$. Then apply the classical $T1$ theorem we see that $S$ is bounded on $L^2$. Hence $T$ is bounded on $L^2$. This completes the proof.

\end{proof}


\begin{thebibliography}{99}

\bibitem{AMV}
E. Airta, H. Martikainen and E. Vuorinen,
\emph{UMD-extensions of Calder\'on--Zygmund operators with mild kernel regularity},
J. Fourier Anal. Appl. \textbf{28} (2022), Paper No.~64.


\bibitem{AY}
P. Auscher and Q.~X. Yang,
\emph{BCR algorithm and the $T(b)$ theorem},
Publ. Mat. \textbf{53} (2009), no.~1, 179--196.


\bibitem{CZ}
A.~P. Calder\'on and A. Zygmund,
\emph{On the existence of certain singular integrals},
Acta Math. \textbf{88} (1952), 85--139.


\bibitem{DJ}
G. David and J.-L. Journ\'e,
\emph{A boundedness criterion for generalized Calder\'on--Zygmund operators},
Ann. of Math. (2) \textbf{120} (1984), no.~2, 371--397.


\bibitem{DYY}
D. Deng, L. Yan and Q. Yang,
\emph{Blocking analysis and $T(1)$ theorem},
Sci. China Ser. A \textbf{41} (1998), no.~8, 801--808.


\bibitem{Figiel}
T. Figiel,
\emph{Singular integral operators: a martingale approach},
in \emph{Geometry of Banach Spaces} (Strobl, 1989),
London Math. Soc. Lecture Note Ser., vol.~158,
Cambridge Univ. Press, Cambridge, 1990, pp.~95--110.


\bibitem{GH}
A. Grau de la Herr\'an and T. Hyt\"onen,
\emph{Dyadic representation and boundedness of nonhomogeneous
Calder\'on--Zygmund operators with mild kernel regularity},
Michigan Math. J. \textbf{67} (2018), no.~4, 757--786.


\bibitem{Lacey}
M.~T. Lacey,
\emph{An elementary proof of the $A_2$ bound},
Israel J. Math. \textbf{217} (2017), no.~1, 181--195.


\bibitem{Lerner}
A.~K. Lerner,
\emph{On pointwise estimates involving sparse operators},
New York J. Math. \textbf{22} (2016), 341--349.


\bibitem{LO}
A.~K. Lerner and S. Ombrosi,
\emph{Some remarks on the pointwise sparse domination},
J. Geom. Anal. \textbf{30} (2020), no.~1, 1011--1027.


\end{thebibliography}
\end{document}